\documentclass[12 pt,a4paper]{amsart}
\usepackage{amssymb}
\usepackage{stmaryrd}
\usepackage{bbm}

\newtheoremstyle{break}{9pt}{9pt}{\itshape}{}{\bfseries}{}{\newline}{}
\newtheoremstyle{rem}{9pt}{9pt}{}{}{\bfseries}{}{\newline}{}

\theoremstyle{break}
\newtheorem*{them}{Theorem}
\newtheorem{lem}{Lemma}
\newtheorem{prop}[lem]{Proposition}

\theoremstyle{rem}
\newtheorem{rem}[lem]{Remark}

\DeclareMathOperator{\Ass}{Ass}
\DeclareMathOperator{\NZD}{NZD}
\DeclareMathOperator{\Var}{Var}
\DeclareMathOperator{\Spec}{Spec}
\DeclareMathOperator{\Min}{Min}
\DeclareMathOperator{\Gr}{Gr}
\DeclareMathOperator{\length}{length}

\newcommand{\Nat}{\mathbbm{N}_0}
\newcommand{\defin}{\mathrel{\mathop:}=}
\newcommand{\Hn}[3]{H^{#1}_{#2}(#3)}
\newcommand{\Hi}[2]{H^{#1}_{R_+}(#2)}
\newcommand{\res}[1]{\!\!\upharpoonright_{#1}}
\newcommand{\finfpot}[2]{#1\llbracket X_1,\ldots,X_{#2}\rrbracket}
\newcommand{\Id}{{\rm Id}}

\begin{document}
		
\title[An Avoidance Principle with an Application to Local Cohomology]{An Avoidance Principle with an\\Application to the Asymptotic Behaviour\\of Graded Local Cohomology}
\author{M. Brodmann}
\author{S. Kurmann}
\author{F. Rohrer}
\address{Institute of Mathematics\\University of Z\"{u}rich\\Winterthurerstrasse 190\\CH-8057 Z\"{u}rich\\Switzerland}
\email{brodmann@math.unizh.ch}
\email{simon.kurmann@math.unizh.ch}
\email{fred.rohrer@math.unizh.ch}
\date{Z\"{u}rich, September 2006}

\begin{abstract}
We present an Avoidance Principle for certain graded rings. As an application we fill a gap in the proof of a result by Brodmann, Rohrer and Sazeedeh about the antipolynomiality of the Hilbert-Samuel multiplicity of the graded components of the local cohomology modules of a finitely generated module over a Noetherian homogeneous ring with two-dimensional local base ring.
\end{abstract}

\maketitle

%%%%%%%%%%%%%%%%%%%%%%%%%%%%%%%%%%%%%%%%%%%%%%%%%%%%%%%%%%%%%%%%%%%%%%%%%%%%%%%%%%%%%%%%%%%%%%%%%%%%%%%%%%

Let $R=\bigoplus_{n\in\Nat}R_n$ be a Noetherian homogeneous ring with local base ring $(R_0,\mathfrak{m}_0)$ and irrelevant ideal $R_+\defin\bigoplus_{n\in\mathbbm{N}}R_n$, let $\mathfrak{q}_0\subseteq R_0$ be an $\mathfrak{m}_0$-primary ideal and let $M=\bigoplus_{n\in\mathbbm{Z}}M_n$ be a finitely generated graded $R$-module. For $i\in\Nat$ we denote by $\Hi{i}{M}$ the $i$-th local cohomology module of $M$ with respect to $R_+$, and for $n\in\mathbbm{Z}$ we denote by $\Hi{i}{M}_n$ its $n$-th graded component. By $e_{\mathfrak{q}_0}(T)$ we denote the Hilbert-Samuel multiplicity of a finitely generated $R_0$-module $T$ with respect to $\mathfrak{q}_0$.

If $A$ is a ring, $\mathfrak{a}\subseteq A$ is an ideal and $T$ is an $A$-module, by $\Gr(\mathfrak{a})\defin\bigoplus_{n\in\Nat}\mathfrak{a}^n/\mathfrak{a}^{n+1}$ and $\Gr(\mathfrak{a},T)\defin\bigoplus_{n\in\Nat}\mathfrak{a}^nT/\mathfrak{a}^{n+1}T$ we denote the graded ring associated with $\mathfrak{a}$ and the graded module associated with $T$ with respect to $\mathfrak{a}$ respectively.

If $(A,\mathfrak{m})$ is a local ring, by an unramified $A$-algebra we mean an $A$-algebra $B$ which is a local ring with maximal ideal $\mathfrak{m}B$.

\medskip

The following Theorem is stated in \cite[Theorem 5.7]{BRS}:
	
\begin{them}\label{Patient}
Let $i\in\Nat$, let $\dim(R_0)=2$ and let $\dim_{R_0}(\Hi{i}{M}_n)\ge 1$ for all $n\ll 0$. Then, there exists a polynomial $Q\in\mathbbm{Q}[X]$ such that $\deg(Q)<i$ and $e_{\mathfrak{q}_0}(\Hi{i}{M}_n)=Q(n)$ for all $n\ll 0$.
\end{them}

But Alas!, there is a gap in the proof of this Theorem in \cite{BRS}. In this note we will fix this. We will not reproduce the entire proof and refer the reader (as we do for unexplained notations) to the article cited above.

The argument lacks in the case where $\dim_{R_0}(\Hi{i}{M}_n)=1$ for all $n\ll 0$ and under the additional assumptions that $\Gamma_{\mathfrak{m}_0R}(M)=0$ and that $R_0$ is complete. For $n\in\mathbbm{Z}$ let $T_n\defin\Hi{i}{M}_n/\Gamma_{\mathfrak{m}_0}(\Hi{i}{M}_n)$. In order to continue the proof there must exist -- possibly after a replacement of $R_0$ by an appropriate ring -- an element $x\in\mathfrak{m}_0\cap\NZD_{R_0}(M)$ which is a parameter for $R_0$ and a non-zerodivisor on $T_n$ for $n\ll 0$ such that for all $n\ll 0$ the following equations hold: $$e_{\mathfrak{q}_0}(\Hi{i}{M}_n)=e_{\mathfrak{q}_0}(T_n)=e_{\mathfrak{q}_0}(T_n/xT_n)=$$$$\length_{R_0}(T_n/xT_n)= \length_{R_0}\bigl(\Hi{i}{M}_n/\bigl(x\Hi{i}{M}_n+\Gamma_{\mathfrak{m}_0}(\Hi{i}{M}_n)\bigr)\bigr).$$ In \cite{BRS} it was not shown that the second equation $e_{\mathfrak{q}_0}(T_n)=e_{\mathfrak{q}_0}(T_n/xT_n)$ holds for the element $x$ given there. We begin with a remark which reminds of some facts on superficial elements (cf. \cite[\S\S 22f.]{na}) in a form suitable to our needs.

\begin{rem}\label{superficial}
Let $A$ be a Noetherian ring, let $\mathfrak{a}\subseteq A$ be an ideal and let $T$ be a finitely generated $A$-module. Let $x\in\mathfrak{a}\cap\NZD_A(T)$ be such that $x+\mathfrak{a}^2\in\Gr(\mathfrak{a})_1=\mathfrak{a}/\mathfrak{a}^2$ is a quasi-non-zerodivisor on $\Gr(\mathfrak{a},T)$, i.e. there is an $m\in\Nat$ such that $x+\mathfrak{a}^2\in\NZD_{\Gr(\mathfrak{a})}(\bigoplus_{k\geq m}\Gr(\mathfrak{a},T)_k)$. Then, there is some $m_0\in\Nat$ such that $$\frac{1}{x}\mathfrak{a}^mT\cap\mathfrak{a}^{m-2}T=(\mathfrak{a}^mT:_{\mathfrak{a}^{m-2}T}x)=\mathfrak{a}^{m-1}T$$ for all $m\geq m_0$. By induction it follows that $\frac{1}{x}\mathfrak{a}^mT\cap\mathfrak{a}^{m_0}T\subseteq\mathfrak{a}^{m-1}T$ for all $m\geq m_0$. By Artin-Rees there is some $m_1\geq m_0$ such that $\mathfrak{a}^{m_1}\frac{1}{x}T\cap T\subseteq\mathfrak{a}^{m_0}T$. So $\frac{1}{x}\mathfrak{a}^mT\cap T\subseteq\mathfrak{a}^{m-1}T$ and hence $(\mathfrak{a}^mT:_Tx)=\mathfrak{a}^{m-1}T$ for all $m\geq m_1$. Therefore, the exact sequences $$0\rightarrow T/(\mathfrak{a}^mT:_Tx)\rightarrow T/\mathfrak{a}^mT\rightarrow(T/xT)/\mathfrak{a}^m(T/xT)\rightarrow 0$$ yield isomorphisms of $A$-modules $(T/xT)/\mathfrak{a}^m(T/xT)\cong(T/\mathfrak{a}^mT)/(T/\mathfrak{a}^{m-1}T)$ for all $m\gg 0$.

If $(A,\mathfrak{m})$ is local, $\mathfrak{a}$ is $\mathfrak{m}$-primary and $\dim_A(T)>0$, it follows for the Hilbert-Samuel polynomials that $P_{\mathfrak{a},T/xT}({\bf x})=P_{\mathfrak{a},T}({\bf x})-P_{\mathfrak{a},T}({\bf x}-1)$ and hence $e_{\mathfrak{a}}(T/xT)=e_{\mathfrak{a}}(T)$.\hfill$\bullet$
\end{rem}

By Remark \ref{superficial}, for the equations $e_{\mathfrak{q}_0}(T_n)=e_{\mathfrak{q}_0}(T_n/xT_n)$ as well as the other conditions to hold it is sufficient for $x$ to meet the following conditions:
\begin{enumerate}
\item[(1)] $x$ is a parameter for $R_0$;
\item[(2)] $\forall\, n\ll 0:x\in\NZD_{R_0}(T_n)$;
\item[(3)] $\forall\, k\in\mathbbm{Z}:x\in\NZD_{R_0}(M_k)$;
\item[(4)] $x\in\mathfrak{q}_0$;
\item[(5)] $\forall\, n\ll 0:x+\mathfrak{q}_0^2\in\Gr(\mathfrak{q}_0)_1$ is a quasi-non-zerodivisor on $\Gr(\mathfrak{q}_0, T_n)$.
\end{enumerate}
We shall indeed look for an element $x$ which satisfies the following stronger conditions:
\begin{enumerate}
\item[(1')] $x\notin\bigcup\Min(R_0)$;
\item[(2')] $\forall\, n\in\mathbbm{Z}:x\notin\bigcup\Ass_{R_0}(T_n)$;
\item[(3')] $\forall\, k\in\mathbbm{Z}:x\notin\bigcup\Ass_{R_0}(M_k)$;
\item[(4')] $x\in\mathfrak{q}_0\setminus\mathfrak{m}_0\mathfrak{q}_0$;
\item[(5')] $\forall\, n\in\mathbbm{Z}:x+\mathfrak{q}_0^2\notin\bigcup\Ass_{\Gr(\mathfrak{q}_0)}\left(\Gr(\mathfrak{q}_0,T_n)\right)\setminus\Var(\Gr(\mathfrak{q}_0)_+)$.
\end{enumerate}	
(For $n\in\mathbbm{Z}$, keep in mind that, as $\Gr(\mathfrak{q}_0)$ is Noetherian and homogeneous and $\Gr(\mathfrak{q}_0,T_n)$ is finitely generated, $x+\mathfrak{q}_0^2$ is a quasi-non-zerodivisor on $\Gr(\mathfrak{q}_0,T_n)$ if and only if it avoids all essential primes associated with $\Gr(\mathfrak{q}_0,T_n)$.)

As $R_0$ is complete, the Countable Prime Avoidance Principle (cf. \cite[Corollary 2.2]{SV}) allows one to find an $x\in\mathfrak{m}_0$ which satisfies conditions (1'), (2') and (3'). This was already used in the proof of \cite[Theorem 5.7]{BRS}. It remains to show that there is an element $x$ which moreover satisfies conditions (4') and (5'). After some reduction steps, this is done by the Avoidance Principle presented below (Proposition \ref{vermeidung}).

\smallskip

To prove our Avoidance Principle we need the following result.

\begin{lem}\label{linalg}
Let $K$ be a field, let $V$ be a vector space over $K$ of finite dimension and let $(V_i)_{i\in I}$ be a family of proper subspaces $V_i\subsetneqq V$ such that $I$ has cardinality strictly less than $K$. Then $\bigcup_{i\in I}V_i\subsetneqq V$.
\end{lem}

\begin{proof}
We argue by induction on $d\defin\dim_K(V)$. The cases $d=0$ and $d=1$ are obviously true. Let $d>1$. Without loss of generality we assume that $(V_i)_{i\in I}$ consists of hyperlanes in $V$. If $(e_j)_{j=1}^d$ is a basis of $V$, the set $\{\langle e_1+ke_2,e_3,\ldots,e_d\rangle_K|k\in K\}$ consists of hyperlanes in $V$ and has the same cardinality as $K$. Hence, there is a hyperplane $W$ in $V$ not occuring in $(V_i)_{i\in I}$. Therefore, $(V_i\cap W)_{i\in I}$ consists of proper subspaces of $W$. The induction hypothesis implies $\bigcup_{i\in I}(V_i\cap W)\subsetneqq W$, and from this we get the claim.
\end{proof}

\begin{prop}\label{vermeidung}
Let $(A,\mathfrak{m})$ be a local ring with infinite residue field $K$ and let $\mathfrak{a}\subseteq A$ be a finitely generated ideal with $\mathfrak{a}\neq 0$. Let $S$ be a set of ideals of $A$ not containing $\mathfrak{a}$ and let $T$ be a set of ideals of $\Gr(\mathfrak{a})$ not containing $\Gr(\mathfrak{a})_+$. Assume that $S\cup T$ has cardinality strictly less than $K$. Then, there is an element $x\in\mathfrak{a}\setminus(\mathfrak{m}\mathfrak{a}\cup\bigcup S)$ such that $x+\mathfrak{a}^2\in\Gr(\mathfrak{a})_1\setminus\bigcup T$.
\end{prop}

\begin{proof}
Consider the familiy $F\defin((\mathfrak{b}\cap\mathfrak{a}+\mathfrak{m}\mathfrak{a})/\mathfrak{m}\mathfrak{a})_{\mathfrak{b}\in S}$ of $K$-subspaces of $\mathfrak{a}/\mathfrak{m}\mathfrak{a}$. Nakayama's Lemma implies $\mathfrak{b}\cap\mathfrak{a}+\mathfrak{m}\mathfrak{a}\subsetneqq\mathfrak{a}$ for each $\mathfrak{b}\in S$. Thus, $F$ consists of proper $K$-subspaces of $\mathfrak{a}/\mathfrak{m}\mathfrak{a}$.

For each $\mathfrak{c}\in T$ there is an ideal $\mathfrak{d}_{\mathfrak{c}}\subseteq A$ such that $\mathfrak{a}^2\subseteq\mathfrak{d}_{\mathfrak{c}}\subseteq\mathfrak{a}$ and $\mathfrak{c}\cap\Gr(\mathfrak{a})_1=\mathfrak{c}\cap(\mathfrak{a}/\mathfrak{a}^2)=\mathfrak{d}_{\mathfrak{c}}/\mathfrak{a}^2$. As $\Gr(\mathfrak{a})_+\not\subseteq\mathfrak{c}$, we have $\mathfrak{d}_{\mathfrak{c}}\subsetneqq\mathfrak{a}$. So, Nakayama's Lemma implies $\mathfrak{d}_{\mathfrak{c}}+\mathfrak{m}\mathfrak{a}\subsetneqq\mathfrak{a}$, and $G\defin((\mathfrak{d}_{\mathfrak{c}}+\mathfrak{m}\mathfrak{a})/\mathfrak{m}\mathfrak{a})_{\mathfrak{c}\in T}$ is a family of proper $K$-subspaces of $\mathfrak{a}/\mathfrak{m}\mathfrak{a}$. As the family $F\amalg G$ has cardinality strictly less than $K$, by Lemma \ref{linalg} we find some $x\in\mathfrak{a}\setminus\mathfrak{m}\mathfrak{a}$ such that $x+\mathfrak{m}\mathfrak{a}\in\mathfrak{a}/\mathfrak{m}\mathfrak{a}$ avoids all members of $F\amalg G$. This element $x$ has the requested properties.
\end{proof}

The reduction will be based on a lifting result for complete local rings. We begin with the following result about the existence of certain flat algebras over a discrete valuation ring (DVR). For its proof we refer the reader to \cite[Proposition 0.6.8.2, Corollaire 0.6.8.3]{EGA}.

\begin{lem}\label{erweiterung}
Let $A$ be a DVR with residue field $K$ and let $L$ be an extension field of $K$. Then, there exists an unramified flat $A$-algebra which is a complete DVR with residue field isomorphic to $L$.
\end{lem}

\begin{prop}\label{lifting}
Let $A$ be a Noetherian complete local ring with residue field $K$ and let $L$ be an extension field of $K$. Then, there is a diagram of rings $$A''\overset{g}\longleftarrow A'\overset{f}\longrightarrow A$$ such that $A'$ and $A''$ are complete regular local rings, $f$ is surjective and $g$ turns $A''$ into an unramified flat $A'$-algebra with residue field isomorphic to $L$.
\end{prop}

\begin{proof}
As $A$ is complete, Cohen's Structure Theorem (cf. \cite[\S 29]{mat}) gives rise to a surjective morphism of rings $f:A'\rightarrow A$, where $A'=\finfpot{B}{t}$ is a ring of formal power series with $t\in\Nat$ and $B$ a field or a complete DVR. In particular, $A'$ is a complete regular local ring, and if $K'$ denotes the residue field of $A'$, we have $K'\cong K$. In order to construct a morphism of the required type $g:A'\rightarrow A''$ we consider the two possible cases separately:

{\it Case 1: $B$ is a field.} Then $B\cong K'\cong K$, so that $L$ is an extension field of $B$. Now, set $A''\defin\finfpot{L}{t}$ and let $g:A'\rightarrow A''$ be the morphism of rings induced by the inclusion $B\subseteq L$. Then, $A''$ is a complete regular local ring and an unramified $A'$-algebra with residue field isomorphic to $L$. Moreover, by the Local Flatness Criterion (cf. \cite[Theorem 22.1]{mat}) $g$ is flat.

{\it Case 2: $B$ is a DVR.} Let $\pi\in B$ be a uniformizer of $B$. As $L$ is an extension field of $B/\pi B$, by Lemma \ref{erweiterung} there is an unramified flat $B$-algebra $C$ which is a complete DVR with residue field isomorphic to $L$. Now, set $A''\defin\finfpot{C}{t}$ and let $g:A'\rightarrow A''$ be the morphism of rings induced by the structure morphism $B\rightarrow C$. Then, $A''$ is a complete regular local ring with residue field isomorphic to $L$ and again in the view of the Local Flatness Criterion $g$ turns $A''$ into an unramified flat $A'$-algebra.
\end{proof}

For the reduction to be applicable we will need the following Lemma. (It does not need any assumptions beyond those in the first paragraph of this article.)

\begin{lem}\label{unverzweigt}
Given a diagram of Noetherian local rings $$R_0''\overset{g_0}\longleftarrow R_0'\overset{f_0}\longrightarrow R_0$$ such that $f_0$ is surjective and $g_0$ turns $R_0''$ into an unramified flat $R_0'$-algebra, for all $i\in\Nat$ and $n\in\mathbbm{Z}$ we have $$e_{\mathfrak{q}_0}(\Hi{i}{M}_n)=e_{f_0^{-1}(\mathfrak{q}_0)R_0''}(\Hn{i}{R_+\otimes_{R_0'}R_0''}{M\otimes_{R_0'}R_0''}_n).$$
\end{lem}

\begin{proof}
As $R$ is Noetherian and homogeneous, there are finitely many elements $l_1,\ldots,l_r$ $\in R_1$ such that $R=R_0[l_1,\ldots,l_r]$. Let $R'\defin R_0'[X_1,\ldots,X_r]$ be the polynomial ring over $R_0'$ furnished with the standard $\Nat$-grading. Then, there is a surjective morphism of graded rings $f:R'\rightarrow R$ with component of degree $0$ equal to $f_0$.

Let $i\in\Nat$ and $n\in\mathbbm{Z}$. As $f$ is surjective, by the graded base ring independence property of local cohomology (cf. \cite[Theorem 13.1.6]{BS}) there is an isomorphism of $R_0'$-modules $\Hi{i}{M}_n\res{R_0'}\cong\Hn{i}{R_+'}{M\res{R'}}_n$ and $f_0^{-1}(\mathfrak{q}_0)$ is $\mathfrak{m}_0$-primary. Now it is easy to see that $$e_{\mathfrak{q}_0}(\Hn{i}{R_+}{M}_n)=e_{f^{-1}(\mathfrak{q}_0)}(\Hn{i}{R_+'}{M\res{R'}}_n).$$ Therefore, it suffices to prove the claim under the assumption that $R_0'=R_0$ and $f_0=\Id_{R_0'}$.

Let $\mathfrak{m}_0''$ denote the maximal ideal of $R_0''$. Note that by unramifiedness $\mathfrak{q}_0R_0''$ is $\mathfrak{m}_0''$-primary. Looking at chains of submodules, by flatness and unramifiedness it is easy to see that $\length_{R_0''}(T\otimes_{R_0}R_0'')=\length_{R_0}(T)$ for an $R_0$-module $T$. Writing $R_+''\defin R_+\otimes_{R_0}R_0''$, by the graded flat base change property of local cohomology (cf. \cite[Theorem 13.1.8]{BS}) we have $$\Hi{i}{M}_n\otimes_{R_0}R_0''\cong\Hn{i}{R_+''}{M\otimes_{R_0}R_0''}_n.$$ Hence, again by flatness, we get an isomorphism $(\Hi{i}{M}_n/\mathfrak{q}_0^{m+1}\Hi{i}{M}_n)\otimes_{R_0}R_0''\cong \Hn{i}{R_+''}{M\otimes_{R_0}R_0''}_n/(f^{-1}_0(\mathfrak{q}_0)^{m+1}R_0'')\Hn{i}{R_+''}{M\otimes_{R_0}R_0''}_n$ for every $m\in\Nat$. The above remark about the lengths and the definition of the Hilbert-Samuel multiplicity now yields the claim.
\end{proof}

%\bigskip

Now, we can use Proposition \ref{lifting} to perform a replacement argument which brings us to a situation in which we have an element $x$ satisfying the conditions (1') to (5'), as the following argument shows.

We choose an uncountable extension field $L$ of $K\defin R_0/\mathfrak{m}_0$, for example the field of fractions of a polynomial ring in uncountably many indeterminates over $R_0/\mathfrak{m}_0$. By Proposition \ref{lifting} there is a diagram of rings $$R_0''\overset{g_0}\longleftarrow R_0'\overset{f_0}\longrightarrow R_0$$ such that $(R_0',\mathfrak{m}_0')$ and $(R_0'',\mathfrak{m}_0'')$ are complete regular local rings with residue fields denoted by $K'$ and $K''$ respectively, $f_0$ is surjective, $g_0$ turns $R_0''$ into an unramified flat $R_0'$-algebra and $K''\cong L$. By Lemma \ref{unverzweigt}, $$e_{\mathfrak{q}_0}(\Hn{i}{R_+}{M}_n)=e_{f_0^{-1}(\mathfrak{q}_0)R_0''}(\Hn{i}{R_+\otimes_{R_0'}R_0''}{M\otimes_{R_0'}R_0''}_n).$$ Therefore, we may replace $R$ and $M$ by $R\otimes_{R_0'}R_0''$ and $M\otimes_{R_0'}R_0''$, respectively, and hence assume that the residue field $K$ of the base ring $R_0$ is uncountable.

Thus, we can use Proposition \ref{vermeidung} as follows. For $n\in\mathbbm{Z}$ we define $$\mathcal{A}\defin\Min(R_0)\cup\bigcup_{n\in\mathbbm{Z}}\Ass_{R_0}(T_n)\cup\bigcup_{k\in\mathbbm{Z}}\Ass_{R_0}(M_k)\subseteq\Spec(R_0)$$ and $$\mathcal{B}\defin\bigcup_{n\in\mathbbm{Z}}\Ass_{\Gr(\mathfrak{q}_0)}(\Gr(\mathfrak{q}_0,T_n))\setminus\Var(\Gr(\mathfrak{q}_0)_+)\subseteq\Spec(\Gr(\mathfrak{q}_0)).$$ As $\Gr(\mathfrak{q}_0)$ is Noetherian and $\Gr(\mathfrak{q}_0,T_n)$ is finitely generated for $n\in\mathbbm{Z}$, both $\mathcal{A}$ and $\mathcal{B}$ are countable sets of ideals of $R_0$ and $\Gr(\mathfrak{q}_0)$ respectively. Furthermore, as $\Gamma_{\mathfrak{m}_0R}(M)=0$, no ideal in $\mathcal{A}$ contains $\mathfrak{q}_0$ and no ideal in $\mathcal{B}$ contains $\Gr(\mathfrak{q}_0)_+$. By Proposition \ref{vermeidung} there is an element $x\in\mathfrak{q}_0\setminus(\mathfrak{m}_0\mathfrak{q}_0\cup\bigcup\mathcal{A})$ such that $x+\mathfrak{q}_0^2\in\Gr(\mathfrak{q}_0)_1\setminus\bigcup\mathcal{B}$. Thus $x$ meets conditions (1') to (5') from above, and Lo!, the gap is filled.

%%%%%%%%%%%%%%%%%%%%%%%%%%%%%%%%%%%%%%%%%%%%%%%%%%%%%%%%%%%%%%%%%%%%%%%%%%%%%%%%%%%%%%%%%%%%%%%%%%%%%%%%%%

\end{document}